\documentclass[12pt]{article}
\pdfoutput=1
\usepackage{amsmath}    % need for subequations
\usepackage{graphicx}   % need for figures
\usepackage{verbatim}   % useful for program listings
\usepackage{color}      % use if color is used in text
\usepackage{subfigure}  % use for side-by-side figures
\usepackage{amssymb}

\usepackage{mathtools}
\DeclarePairedDelimiter\ceil{\lceil}{\rceil}

\newtheorem{theorem}{Theorem}[section]

\newtheorem{corollary}[theorem]{Corollary}

\newtheorem{lemma}[theorem]{Lemma}

\newtheorem{property}[theorem]{Property}

\def\cp{\,\square\,}

\begin{document}
	
\title{On the isometric path partition problem}
	
\author{
	Paul Manuel
}
{\Large This is a part of a proposed research project of Kuwait University, Kuwait. The author will appreciate to receive your comments and constructive criticism on this initial draft. The author's contact address is pauldmanuel@gmail.com.}
\date{}

\maketitle
\vspace{-0.8 cm}
\begin{center}
	Department of Information Science, College of Computing Science and Engineering, Kuwait University, Kuwait \\
	{\tt pauldmanuel@gmail.com}\\
	\medskip
		
\end{center}

\begin{abstract}
The isometric path cover (partition) problem of a graph is to find a minimum set of isometric paths which cover (partition) the vertex set of the graph. The isometric path cover (partition) number of a graph is the cardinality a minimum  isometric path cover (partition). We prove that the isometric path partition problem and the isometric $k$-path partition problem for $k\geq 3$ are NP-complete on general graphs. Fisher and Fitzpatrick \cite{FiFi01} have shown that the isometric path cover number of $(r\times r)$-dimensional grid is $\lceil 2r/3\rceil$. We show that the isometric path cover (partition) number of $(r\times s)$-dimensional grid is $s$ when $r \geq s(s-1)$. We  establish that the isometric path cover (partition) number of $(r\times r)$-dimensional torus is $r$ when $r$ is even and is either $r$ or $r+1$ when $r$ is odd. Then, we demonstrate that the isometric path cover (partition) number of  an $r$-dimensional Benes network is $2^r$. In addition, we provide partial solutions for the isometric path cover (partition) problems for cylinder and multi-dimensional grids. 
\end{abstract}


\medskip
\noindent{\bf AMS Subj.\ Class.: 05C12, 05C70, 68Q17}

%%%%%%%%%%%%%%%%%%%%%%%%%%%
\section{Introduction}
%%%%%%%%%%%%%%%%%%%%%%%%%%%
An undirected connected graph is represented by $G(V,E)$ where $V$ is the vertex set and $E$ is the edge set. A \emph{path} means a simple path with distinct vertices. A path $P$ is an \emph{induced path} in $G$ if the subgraph induced by the vertices of $P$ is a path. A path between two vertices is \emph{an isometric path} if it induces the shortest distance between the two points. Let us recall that \emph{isometric path} and geodesic are other names for shortest path. While an \emph{isometric path cover} is a set of isometric paths which cover the vertex set $V$, an \emph{isometric path partition} is a set of isometric paths which partition $V$.  The \emph{isometric path cover  number}  which is denoted by ${\rm ip}_{c}(G)$  is the cardinality of a minimum isometric path cover. The isometric path partition number ${\rm ip}_{p}(G)$ is defined accordingly. While the \emph{isometric path cover  problem} is to find a minimum isometric path cover, the \emph{isometric path partition  problem} is to find a minimum isometric path partition. A \emph{diametral isometric path} of a graph is an isometric path whose length is equal to the diameter of the graph. As the path cover (partition) problem is to find a minimum path cover (partition), the  induced path cover (partition) problem is to find a minimum induced path cover (partition). Let ${\rm diam}(G)$ denote the diameter of graph $G$. 

Last few decades, the theory of isometric paths has been studied extensively. Aggarwal et al. \cite{Aggarwal00} have illustrated the application of the isometric path cover problem in the design of VLSI layouts. The theory of isometric path problems is the backbone in the design of efficient algorithms in 
transport networks \cite{OrLl14}, computer networks \cite{LaFe16, PiMe04}, parallel architectures \cite{Xu2013}, social networks \cite{BaTh10, GoLi16}, VLSI layout design \cite{Aggarwal00}, wireless sensor networks \cite{CoRi16}, multimedia networks \cite{ClCr03} and in other networks such as GIS networks \cite{ZhLu10}, large network systems \cite{AkIw13} and stochastic networks \cite{PeSh07}. 

Since the Hamiltonian path problem is NP-complete \cite{GaJo79}, the path cover problem and the path partition problem are NP-complete. Since the Hamiltonian induced path problem is NP-complete \cite{Char94, GaJo79}, the induced path cover problem and the induced path partition problem are NP-complete. However, the complexity status of the isometric path cover problem and the isometric path partition problem are unknown \cite{Manuel18}. This fact has been  highlighted and emphasized recently \cite{Manuel18, MaKlAn18}.  In this paper, we settle the long-standing open problem \cite{Manuel18} by proving that the isometric path partition problem is NP-complete on general graphs. 

The isometric path cover number has been computed for tress, cycles, complete bipartite graphs, the Cartesian product of paths (including hypercubes) under some restricted cases \cite{FiFi01, Fitz97, Fitz99, FiNo01}. Fisher and Fitzpatrick \cite{FiFi01} have derived a lower bound that ${\rm ip}_{c}(G) \geq \ceil*{\frac{|V|}{{\rm diam}(G)+1}}$ and have shown that the isometric path cover number of $(r \times r)$ grid is $\ceil*{{2r}/{3}}$. Fitzpatrick et al. \cite{FiNo01} have shown that the isometric path cover number of hypercube $Q_r$ is at least $2r/(r+1)$. In addition, they have also shown that ${\rm ip}_{c}(Q_r) = 2^{r-log_{2}(r+1)}$ when $r+1$ is a power of 2. Pan and Chang have given a linear-time algorithm to solve the isometric path cover problem on block graphs \cite{PaCh05a}, complete $r$-partite graphs and Cartesian products of 2 or 3 complete graphs \cite{PaCh06}. There is no literature on the isometric path partition problem \cite{Manuel18}. The readers are suggested to read the survey paper by Manuel \cite{Manuel18} for detailed information. 

In section \ref{sec:NP_complete}, we show that the isometric path partition problem is NP-complete. We also show that the isometric $k$-path partition problem is NP-complete on general graphs for $k \geq 3$. In section \ref{sec:iso_path_part_grids}, we compute the exact values of the isometric path cover and isometric path partition number for cylinders and multi-dimensional grids under certain conditions.  In sections \ref{sec:torus} and \ref{sec:benes}, we  also derive the exact values of the isometric path cover number and isometric path partition number for square torus and Benes networks.
%%%%%%%%%%%%%%%%%%%%%%%%%%%%%%
\section{The isometric path partition problem is NP-complete for general graphs}
\label{sec:NP_complete}
\begin{figure}[ht!]
	\begin{center}
		\scalebox{0.90}{\includegraphics{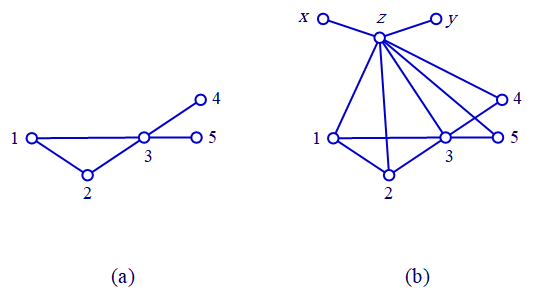}}
	\end{center}
	\caption{(a) $G = (V,E)$ \quad and \quad (b) $G'=(V',E')$}
	\label{fig:FIsoPathPartNP}
\end{figure}

%%%%%%%%%%%%%%%%%%%%%%%%%%%%%%%
Given a graph $G(V,E)$, a \emph{$k$-path} is a path having at most $k$ vertices. A set $S$ of $k$-paths is a \emph{$k$-path partition} if each vertex of $V$ belongs to exactly one member of $S$. The \emph{$k$-path partition problem} is to find a $k$-path partition of minimum cardinality in $G$.
\begin{theorem}[\cite{MoTo07}]
	\label{TPathPartNP1}
	The $3$-path partition problem is {\rm NP}-complete on bipartite graphs. 
\end{theorem}
Now we will prove that the isometric path partition problem is NP-complete on general graphs. As a first step, we provide a polynomial reduction from the $3$-path partition problem on bipartite graphs to the isometric path partition problem on general graphs. The key fact in a bipartite graph is that each $3$-path in a bipartite graph is an isometric path. Given a graph $G(V,E)$ where $V = \{1, 2 \ldots n\}$, the reduced graph is denoted by $G'(V',E')$. The vertex set $V'$ is $V \cup \{x, y, z\}$.  The edge set $E'$ is $E \cup \{xz, zy\} \cup \{iz \,/\, i \in V\}$. See Figure \ref{fig:FIsoPathPartNP}.

Next, we will identify some basic structural properties of $G'(V',E')$. 
\begin{property}
	\label{pro:Pro1}
	If $G$ is a bipartite graph, then ${\rm diam}(G')$ = $2$. For $i, j, k \in V$, a diametral isometric path in $G'$ is either $xzy$, $izx$, $izy$, $ikj$ or $izj$ when $ij$ is not an edge.    
\end{property}

\begin{property}
	\label{pro:Pro2}
	A bipartite graph $G$ has a $3$-path partition $S$ of cardinality $k$ iff  $G'$ has an isometric path partition $S'$ of cardinality $k+1$.   
\end{property}
Applying Property \ref{pro:Pro1} and \ref{pro:Pro2}, we state one of the main results of this paper.
\begin{theorem}
	The isometric path partition problem is NP-complete on general graphs. 
\end{theorem}
The isometric $k$-path partition problem is a generalization of the isometric path partition problem. Using the same logic, one can also prove that
\begin{theorem}
	The isometric $k$-path partition problem is NP-complete on general graphs for $k \geq 3$. 
\end{theorem}
%%%%%%%%%%%%%%%%%%%%%%%%%%%%%%%%
\section{Isometric path partition versus isometric path cover}
\label{sec:cover_vs_partition}
%%%%%%%%%%%%%%%%%%%%%%%%%%%%%%%%%
Our point of discussion in this section is to emphasize that the isometric path cover problem and the isometric path partition problem are two different combinatorial problems.   From the perspective of computational complexity, let us see how these two problems differ even on simple architectures such as trees and grids. For the star graph $K_{1,\ell}$, while ${\rm ip}_{c}(K_{1,\ell})$ = $\ceil*{\ell/2}$, ${\rm ip}_{p}(K_{1,\ell})$ = $\ell-1$. Pan and Chang \cite{PaCh05a} have proved that the isometric path cover number of trees with $\ell$ leaves is $\ceil*{\ell/2}$. For the isometric path partition number, this is not true even in complete binary trees. 
In fact, it is a challenge to find the isometric path partition number for trees and the isometric path partition number for trees is unknown \cite{Manuel18}. 

Fisher and Fitzpatrick \cite{FiFi01} have shown that the isometric path cover number of grid $\Xi(r, r)$ is $\ceil*{2r/3}$. However, the isometric path partition number of grid $\Xi(r, r)$  remains an open problem \cite{Manuel18}. It seems that the isometric path partition number of grid $\Xi(r, r)$  is $r$. However, it requires a mathematical proof which seems to be a challenge.  

In the same way, one can apply greedy algorithm or brute-force algorithm to compute the isometric path cover number for some graphs such as interval graph or circular arc graph by starting from its diametral isometric path. On interval graphs , our experiments demonstrate that  greedy algorithm provides sharp result for the isometric path cover number but significantly deviates for the isometric path partition number.

There are some features which are common to both the problems. One straightforward lower bound that is common to both the problems is as follows:
\begin{theorem}[\cite{FiFi01}]
	\label{thm:trivual_bnd}
	If ${\rm diam}(G)$ denotes the diameter of a graph $G$, then ${\rm ip}_{p}(G) \geq {\rm ip}_{c}(G) \geq \ceil*{ 
	\frac{|V(G)|}{{\rm diam}(G)+1}}$.
\end{theorem}
Though the bound in Theorem \ref{thm:trivual_bnd} seems to be trivial, it is very effective and useful. This paper truly exploits this lower bound to prove that the isometric path cover number and the isometric path partition number are equal for some networks such as $(r \times r)$-dimensional torus and Benes networks. 

Monnot and Toulouse \cite{MoTo07} have studied an NP-complete problem that is whether a graph on $nk$ vertices can be partitioned into $n$ paths of length $k$. An isometric path version of this problem is to decide if a graph $G(V,E)$ can be partitioned into  diametral isometric paths. Any fixed interconnection network which possesses this feature is considered as a ``good'' architecture \cite{Leighton1992, Xu2013} because it is easy to design and implement efficient data communication, message broadcasting and other routing   algorithms. In the following sections, we point out that a few fixed interconnection networks such as torus and Benes networks inherit this nice feature.
%%%%%%%%%%%%%%%%%%%%%%%%%%%%%%%%
\section{The isometric path cover / partition problem on Cartesian product $P_r\cp G$}
\label{sec:iso_path_part_grids}
%%%%%%%%%%%%%%%%%%%%%%%%%%%%%%%%%
Given a graph $G$, let ${\rm diam}(G)$ denote the diameter of $G$. In this section,  we consider the Cartesian product $P_r\cp G$ where $G$ is any graph and $P_r$ is a path graph on $r$ vertices.  The vertex set $V(P_r)$ of $P_r$ is $\{1, 2 \ldots r\}$ and the vertex set $V(P_r\cp G)$ of $P_r\cp G$ is $\{(j,v) \,/\, j \in V(P_r)$ and $ v\in V(G)\}$. For each $j\in V(P_r)$, the subgraph induced by the vertices $\{(j,v) \,/\, v \in V(G)\}$ of $P_r\cp G$ is denoted by $G^j$. In other words, there are $r$ copies of $G$ in $P_r\cp G$ which are represented by $G^1, G^2 \ldots G^r$ respectively. An edge of $G^j$ in $P_r\cp G$ is called $G^j$-edge. A $G$-edge in $P_r\cp G$ is a $G^j$-edge for some $j = 1, 2 \dots r$. 

\begin{lemma}[\cite{ImKl08}]
	\label{PrxG_Edges}
	Given a graph $G$ and a path graph $P_r$, an isometric path of the Cartesian product $P_r\cp G$ can have a maximum of ${\rm diam}(G)$ number of $G$-edges.
\end{lemma}
	
\begin{lemma}
	\label{ip_c_PnxG_UB}
	Given a graph $G$ and a path graph $P_r$, ${\rm ip}_{c}(P_r\cp G)$ $\leq$ ${\rm ip}_{p}(P_r\cp G)$ $\leq$ $|V(G)|$.
\end{lemma}

The following lemma is the key to derive a lower bound on ${\rm ip}_{c}(P_r\cp G))$ and ${\rm ip}_{p}(P_r\cp G))$:
\begin{lemma}
	\label{ip_c_PnxG_LB1}
		Given a graph $G$ and a path graph $P_r$, let $S$ be an isometric path cover (partition) of Cartesian product $P_r\cp G$. If there exists a $G^{j_0}$ of $P_r\cp G$ for some $j_0$, $1 \leq j_0 \leq r$, such that no isometric path of $S$ contains any $G^{j_0}$-edge, then $|S| \geq |V(G)|$.
\end{lemma}

\begin{lemma}
	\label{ip_c_PnxG_LB2}
	Given a graph $G$ and a path graph $P_r$, ${\rm ip}_{p}(P_r\cp G)$ $\geq$ ${\rm ip}_{c}(P_r\cp G)$ $\geq |V(G)|$ when $r \geq {\rm diam}(G)|V(G)|$.
\end{lemma}

Following Lemma \ref{ip_c_PnxG_LB2} and Lemma \ref{ip_c_PnxG_UB}, we state that
\begin{theorem}
	\label{thm:iso_path_cov_Num}
		Given a graph $G$ and a path graph $P_r$, ${\rm ip}_{p}(P_r\cp G))$ = ${\rm ip}_{c}(P_r\cp G))$ = $|V(G)|$ when $r \geq {\rm diam}(G)|V(G)|$.
\end{theorem}

%%%%%%%%%%%%%%%%%%%%%%%%%%%%%%%%
\subsection{The isometric path cover / partition problem on grids and cylinders}
%%%%%%%%%%%%%%%%%%%%%%%%%%%%%%%%
In this section, we consider multi-dimensional grids  $\Xi(d_1, d_2 \ldots d_r)$ which are the Cartesian product of paths $P_{d_1}, P_{d_2} \ldots P_{d_r}$.  Fisher and Fitzpatrick \cite{FiFi01} have shown that the isometric path cover number of $\Xi(r, r)$ grid is $\lceil 2r/3\rceil$. Fitzpatrick et al. \cite{FiNo01} have shown that the lower bound of the isometric path cover number of hypercube $Q_r$ is $2r/(r+1)$. In addition, they have also shown that ${\rm ip}_{c}(Q_r) = 2^{r-log_{2}(r+1)}$ when $r+1$ is a power of 2.
Though there are some results available for the isometric path cover number on grids, there are no literature for the isometric path partition problem on  multi-dimensional grids including 2-dimensional grids and cylinders \cite{Manuel18}.

\begin{theorem}
	\label{thm:ip_c_nulti_grids}
	Given an $r$-dimensional grid $\Xi(d_1, d_2 \ldots d_r)$, ${\rm ip}_{p}(\Xi(d_1, d_2 \ldots d_r))$ = ${\rm ip}_{c}(\Xi(d_1, d_2 \ldots d_r))$ = $d_2 d_3 \ldots d_r$ when $d_1 \geq ((d_2-1)+(d_3-1) \ldots (d_r-1)) (d_2d_3 \ldots d_r)$.
\end{theorem}

\begin{corollary}
	The isometric path cover (partition) number of $(r\times s)$-dimensional grid is $s$ when $r \geq s(s-1)$.
\end{corollary}
\begin{theorem}
	\label{thm:ip_c_cyl}
	Given a cylinder $P_r \cp C_s$, ${\rm ip}_{p}(P_r \cp C_s)$ = ${\rm ip}_{c}(P_r \cp C_s)$ = $s$ when $r \geq \lfloor s/2\rfloor s$.
\end{theorem}

\section{The isometric path partition problem on torus}
\label{sec:torus} 
\begin{figure}[ht!]
	\begin{center}
		\scalebox{0.9}{\includegraphics{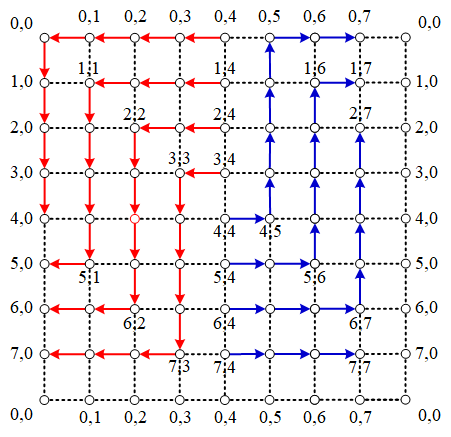}}
	\end{center}
	\caption{An $(8 \times 8)$-dimensional torus}
	\label{fig:square-torus-even}
\end{figure}

\begin{figure}[ht!]
	\begin{center}
		\scalebox{0.9}{\includegraphics{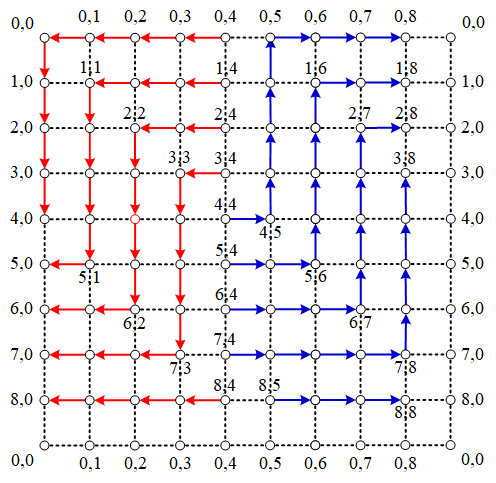}}
	\end{center}
	\caption{A $(9 \times 9)$-dimensional torus}
	\label{fig:square-torus-odd}
\end{figure}

%%%%%%%%%%%%%%%%%%%%%%%%%%%%%%%
In this section, we study the exact value of the isometric path partition number of $(r \times r)$-dimensional torus.  To our knowledge, there is no literature on the isometric path partition problem on torus. In this section, we will show that the isometric path cover (partition) number of $(r\times r)$-dimensional torus $G$ is $r$ when $r$ is even and is either $r$ or $r+1$ when $r$ is odd. An $(8 \times 8)$-dimensional torus is given in  Figure \ref{fig:square-torus-even} and a $(9\times 9)$-dimensional torus is given in Figure \ref{fig:square-torus-odd}.

\begin{theorem}
	\label{thm:torus}
	The isometric path cover $($partition$)$ number of $(r\times r)$-dimensional torus is $r$  when $r$ is even and is either $r$ or $r+1$ when $r$ is odd.
\end{theorem}
 
%%%%%%%%%%%%%%%%%%%%%%%%%%%%%%
\section{The isometric path partition problem on Benes networks}
\label{sec:benes}
%%%%%%%%%%%%%%%%%%%%%%%%%%%%%%%
Let $\mathbb{Z}_k = \{0, 1 \ldots k-1\}$ and $\mathbb{Z}_{2}^{k} =\{x_0x_1 \ldots x_{k-1} / \, x_i = 0$ or $1\}$. When we say $i \in \mathbb{Z}_{k}$, it means $i\mod k$. The vertex set of an $r$-dimensional Butterfly $BF(r)$ is $\{\langle w,i\rangle$ / $w \in \mathbb{Z}_{2}^{r}$ and $i \in \mathbb{Z}^{r+1}\}$. Two vertices $\langle w,i\rangle $ and $\langle w',i'\rangle $ of $BF(r)$ are linked by an edge if $i'=i+1$ and either $w$ = $w'$ or $w$ and $w'$ differ only in the bit in position $i$. 
\begin{figure}[ht!]
	\begin{center}
		\scalebox{0.43}{\includegraphics{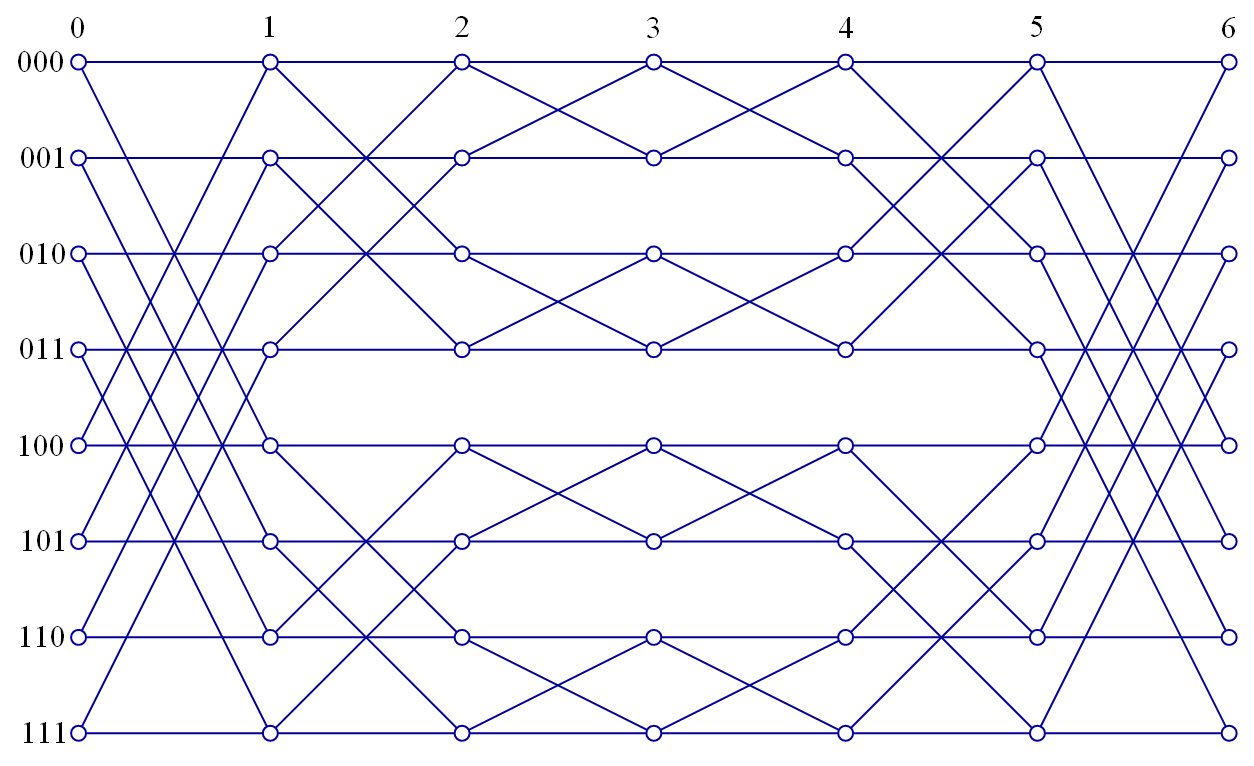}}
	\end{center}
	\caption{A 3-dimensional Benes network }
	\label{fig:Benes-3-dim-nor}
\end{figure}
An $r$-dimensional Benes  network $BN(r)$ is called back-to-back butterflies and is obtained by merging the $r$-level vertices of two Butterfly networks $BF(r)$. The vertex set of  $BN(r)$ is $\{\langle w,i\rangle$ / $w \in \mathbb{Z}_{2}^{r}$ and $i \in \mathbb{Z}^{2r+1}\}$. Thus, the degree of 0-level vertices of $BN(r)$ remains $2$. Figure \ref{fig:Benes-3-dim-nor} displays $BN(3)$. 
We apply the lower bound of Theorem \ref{thm:trivual_bnd} to find the exact value  for ${\rm ip}_{p}(BN(r))$ and ${\rm ip}_{c}(BN(r))$ of Benes networks.

\begin{theorem}
	\label{thm:benes}
	The isometric path cover (partition) number of  an $r$-dimensional Benes network is $2^r$.
\end{theorem}
%%%%%%%%%%%%%%%%%%%%%%%%%%%%%
\section*{Acknowledgment}
%%%%%%%%%%%%%%%%%%%%%%%%%%%%%

This work is a part of research project proposal which is submitted to Kuwait University, Kuwait.


\begin{thebibliography}{10}
	
	\bibitem{Aggarwal00}
	Alok Aggarwal, Jon~M Kleinbergy, and David~P Williamson, \emph{Node-disjoint
		paths on the mesh and a new trade-off in {VLSI} layout}, {SIAM} {J}.
	{C}omput. \textbf{29} (2000), no.~4, 1321--1333.
	
	\bibitem{AkIw13}
	Takuya Akiba, Yoichi Iwata, and Yuichi Yoshida, \emph{Fast exact shortest-path
		distance queries on large networks by pruned landmark labeling}, {SIGMOD} '13
	Proceedings of the 2013 {ACM SIGMOD} International Conference on Management
	of Data, New York, USA, {ACM} Digital Library, 2013, pp.~349--360.
	
	\bibitem{BaTh10}
	John Baras, George Theodorakopoulos, and Jean Walrand, \emph{Path problems in
		networks}, Synthesis Lectures on Communication Networks, Morgan \& Claypool,
	2010.
	
	\bibitem{Char94}
	Gary Chartrand, Joseph McCanna, Naveed Sherwani, Moazzem Hossain, and Jahangir
	Hashmi, \emph{The induced path number of bipartite graphs}, Ars Combin.
	\textbf{37} (1994), 191--208.
	
	\bibitem{ClCr03}
	Jo|v{a}o C.~N. Clímaco, Jos\'{e} M.~F. Craveirinha, and Marta M.~B. Pascoal,
	\emph{A bicriterion approach for routing problems in multimedia networks},
	Network, An International Journal \textbf{41} (2003), no.~4, 206--220.
	
	\bibitem{CoRi16}
	Juan Cota-Ruiz, Pablo Rivas-Perea, Ernesto Sifuentes, and Rafael
	Gonzalez-Landaeta, \emph{A recursive shortest path routing algorithm with
		application for wireless sensor network localization}, {IEEE} Sensors Journal
	\textbf{16} (2016), no.~11, 4631--4637.
	
	\bibitem{FiFi01}
	David~C. Fisher and Shannon~L. Fitzpatrick, \emph{The isometric number of a
		graph}, Journal of Combinatorial Mathematics and Combinatorial Computing
	\textbf{38} (2001), 97--110.
	
	\bibitem{Fitz97}
	Shannon~L. Fitzpatrick, {PhD} thesis, Department of Mathematics, Dalhousie
	University, Nova Scotia, Canada, 1997.
	
	\bibitem{Fitz99}
	Shannon~L. Fitzpatrick, \emph{The isometric path number of the cartesian product of paths},
	Congr. Numer. \textbf{137} (1999), 109--119.
	
	\bibitem{FiNo01}
	Shannon~L. Fitzpatrick, Richard~J. Nowakowski, Derek~A. Holton, and Ian Caines,
	\emph{Covering hypercubes by isometric paths}, Discrete Mathematics
	\textbf{240} (2001), 253--260.
	
	\bibitem{GaJo79}
	Michael Garey and David~S. Johnson, \emph{Computers and {I}ntractability: A
		guide to the theory of {NP}-completeness}, Freeman, New York, 1979.
	
	\bibitem{GoLi16}
	Maoguo Gong, Guanjun Li, Zhao Wang, Lijia Ma, and Dayong Tian, \emph{An
		efficient shortest path approach for social networks based on community
		structure}, CAAI Transactions on Intelligence Technology \textbf{1} (2016),
	no.~1, 114--123.
	
	\bibitem{ImKl08}
	W.~Imrich, S.~Klav\v{z}ar, and D.~F. Rall, \emph{{Topics in Graph Theory}:
		{G}raphs and their {C}artesian product}, A K Peters, Ltd., Wellesley, MA,
	2008.
	
	\bibitem{LaFe16}
	Mohamed~Lamine Lamali, Nasreddine Fergani, Johanne Cohen, and Helia Pouyllau,
	\emph{Path computation in multi-layer networks: Complexity and algorithms},
	{IEEE INFOCOM} 2016 - {IEEE} {I}nternational Conference on Computer
	Communications, 10-15 April 2016, {S}an {F}rancisco, {CA, USA}, {IEEE}
	Computer Society, 2016.
	
	\bibitem{Leighton1992}
	F.~Thompson Leighton, \emph{Introduction to parallel algorithms and
		architectures: arrays{\textperiodcentered} trees{\textperiodcentered}
		hypercubes}, Morgan Kaufmann Publishers, 1992.
	
	\bibitem{Manuel18}
	Paul Manuel, \emph{Revisiting path-type covering and partitioning problems},
	Manuscript (2018).
	
	\bibitem{MaKlAn18}
	Paul Manuel, Sandi Klavzar, Antony Xavier, Andrew Arokiaraj, and Elizabeth
	Thomas, \emph{Strong geodetic problem in networks}, Discussiones Mathematicae
	Graph Theory (2018), 1--15.
	
	\bibitem{Manuel2008}
	Paul~D. Manuel, Mostafa~I. Abd-El-Barr, Indra Rajasingh, and Bharati Rajan,
	\emph{An efficient representation of {B}enes networks and its applications},
	Journal of Discrete Algorithms \textbf{6} (2008), no.~1, 11--19.
	
	\bibitem{MoTo07}
	J\'{e}r\H{o}me Monnot and Sophie Toulouse, \emph{The path partition problem and
		related problems in bipartite graphs}, Journal Operations Research Letters
	\textbf{35} (2007), no.~5, 677--684.
	
	\bibitem{OrLl14}
	Hector Ortega-Arranz, Diego~R. Llanos, and Arturo Gonzalez-Escribano, \emph{The
		shortest-path problem: Analysis and comparison of methods}, Synthesis
	Lectures on Theoretical Computer Science, Morgan \& Claypool, 2014.
	
	\bibitem{PaCh05a}
	Jun-Jie Pan and Gerard~J. Chang, \emph{Isometric path numbers of block graphs},
	Information Processing Letters \textbf{93} (2005), 99--102.
	
	\bibitem{PaCh06}
	Jun-Jie Pan and Gerard~J. Chang, \emph{Isometric path numbers of graphs}, Discrete Mathematics
	\textbf{306} (2006), no.~17, 2091--2096.
	
	\bibitem{PeSh07}
	S.~K. Peer and Dinesh~K. Sharma, \emph{Finding the shortest path in stochastic
		networks}, Computers \& Mathematics with Applications \textbf{53} (2007),
	no.~5, 729--740.
	
	\bibitem{PiMe04}
	Michal Pioro and Deepankar Medhi, \emph{Routing, flow, and capacity design in
		communication and computer networks}, Networkingq, Morgan Kaufmann, 2004.
	
	\bibitem{Xu2013}
	Junming Xu, \emph{Topological structure and analysis of interconnection
		networks}, Network theory and applications, vol.~7, Springer Science \&
	Business Media, 2013.
	
	\bibitem{ZhLu10}
	Jiyi Zhang, Wen Luo, Linwang Yuan, and Weichang Mei, \emph{Shortest path
		algorithm in {GIS} network analysis based on clifford algebra}, 2010 2nd
	International Conference on Future Computer and Communication, 21-24, May
	2010, Wuha, China, {IEEE} Computer Society, 2010.
	
\end{thebibliography}
\end{document}